\documentclass[12pt]{article}
\usepackage{pstricks,pst-node,pst-tree}
\usepackage{amsmath,amsthm}
\usepackage{amssymb}
\usepackage{color,xcolor}
\usepackage{xspace}
\usepackage[colorlinks=true,
linkcolor=green,
filecolor=brown,
citecolor=green]{hyperref}

\def\v{\vert}
\def\n{n}
\def\a{\ensuremath{\mathcal A}\xspace}
\def\p{\ensuremath{\mathcal P}\xspace}

\def\mbf#1{\mathchoice{\hbox{\boldmath $\displaystyle #1$}}
        {\hbox{\boldmath $\textstyle #1$}}
        {\hbox{\boldmath $\scriptstyle #1$}}
        {\hbox{\boldmath $\scriptscriptstyle #1$}}} 

\begin{document} 
\newtheorem*{theorem*}{Theorem}
\newtheorem{theorem}{Theorem}
\newtheorem{defn}[theorem]{Definition}
\newtheorem{lemma}[theorem]{Lemma}
\newtheorem{prop}[theorem]{Proposition}
\newtheorem{cor}[theorem]{Corollary}
\begin{center}
{\Large
A Combinatorial Interpretation for Sequence A345973 in OEIS}    \\  
\vspace{5mm}
{\large David Callan}  \\[2mm]
August 10, 2021
\end{center}

\begin{abstract}
We give a combinatorial interpretation in terms of bicolored ordered trees for the sequence $(a_n)_{n\ge 1}=(1, 1, 1, 2, 3, 6, 10, 20, 36, 73,\dots )$, A345973 in OEIS, whose generating function satisfies the defining identity $\sum_{n\ge 1}a_n x^n = x+\frac{x^2}{\prod_{n\ge 1}(1-a_nx^n)}$.
\end{abstract}

Our purpose is to give a combinatorial interpretation of the sequence  $(a_n)_{n\ge 1}=(1, 1, 1, 2, 3, 6,$ \\
$ 10, 20, 36, 73,\dots )$, \htmladdnormallink{A345973}{http://oeis.org/A345973} in OEIS \cite{oeis}, with defining identity $\sum_{n\ge 1}a_n x^n = x+\frac{x^2}{\prod_{n\ge 1}(1-a_nx^n)}$. Our counting sequence begins (0,1,1,2,\dots ), differing only in the first term.
All our trees have edges colored black or gray and size measured by number of black edges.
Consider the set \a of such bicolored trees, unlabeled but rooted with ordered edges, that satisfy the following 4 conditions: \\
For each nonleaf vertex $v$,
\begin{enumerate}
\vspace*{-4mm}
\item $v$ has at least two black child edges,
\item the gray child edges of $v$ (if any) are right justified, that is, they lie to the right of all black child edges of $v$,
\item the subtrees of $v$ taken left to right have weakly decreasing sizes, 
\end{enumerate}
\vspace*{-4mm}
and for each leaf,
\begin{enumerate}
\setcounter{enumi}{3}
\vspace*{-3mm}
\item its parent edge (if it has one) is black.
\end{enumerate}
\vspace*{-2mm}                    

The  tree on the left in Figure 1 below is certainly not in \a because it violates each of the first three conditions with $v$ taken as the root and it also violates the fourth condition due to the leaf terminating the gray edge,
while the other two trees are in $\a$. Let $\a_n$ denote the set of trees in \a of size $n$.
\vspace*{-5mm}
\begin{center}
\pstree[nodesep=0pt,levelsep=5ex]{\psset{linewidth=.06cm}\Tdot }{\psset{linecolor=gray}\psset{linewidth=.04cm}
        \Tdot \psset{linecolor=black}\psset{linewidth=.06cm}        
        \pstree{ \Tdot }{
                     \Tdot
                     \Tdot
                              }       
                    }
\hspace*{20mm}
\pstree[nodesep=0pt,levelsep=5ex]{\psset{linewidth=.06cm}\Tdot }{\psset{linewidth=.06cm}          
        \pstree{ \Tdot }{
                     \Tdot
                     \Tdot
                              }
         \Tdot  
         \Tdot      
                    }
\hspace*{20mm}
\pstree[nodesep=0pt,levelsep=5ex]{\psset{linewidth=.06cm}\Tdot }{
        \psset{linewidth=.06cm}        
        \pstree{ \Tdot }{
                   \pstree{ \Tdot }{
                     \Tdot
                     \Tdot}
                     \Tdot
                     \Tdot}
        \pstree{ \Tdot }{
                     \Tdot
                     \Tdot}
        \psset{linecolor=gray}\psset{linewidth=.04cm}\pstree{ \Tdot }
                     {\psset{linecolor=black} \psset{linewidth=.06cm}
                     \Tdot
                     \Tdot}
        \psset{linecolor=gray}\psset{linewidth=.04cm}\pstree{ \Tdot }
                    {\psset{linecolor=black} \psset{linewidth=.06cm}
                     \Tdot
                     \Tdot}                                              
                    }
\end{center}
\vspace*{-2mm}
\hspace*{10mm} not in \a \hspace*{19mm} a tree in $\a_{5}$ \hspace*{21mm} a tree in $\a_{13}$ with two gray edges\\[4mm]
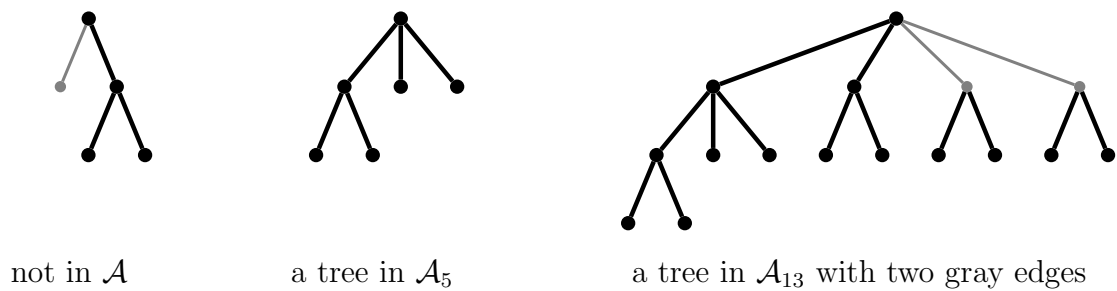
\centerline{Figure 1} 
\newpage
Clearly, $\a_0$ consists of the singleton tree (a root with no edges) and $\a_1$ is empty due to condition 1.
So $\v \a_0 \v=1$ and $\v \a_1 \v=0$.
Set $a_n =\v \a_n \v$ for $n\ge 2$ (leaving $a_1$ undefined for now).
\begin{theorem*}
For $n\ge 2,\  a_n =  A345973(n)$.
\end{theorem*}
\vspace*{-5mm}
\begin{proof}
Suppose $T \in \a_n$ with $n\ge2$.
Say there are $k$ black edges from the root. 
They contribute $k$ to the size of $T$ and $k\ge 2$ by condition 1. The subtrees of the root are all in $\a$ since they inherit the 4 defining conditions.
If there are no gray edges from the root, then either all subtrees of the root are singletons or the singleton subtrees (if any) are right justified by condition 3 and
the non-singleton subtrees, taken left to right, have sizes (weakly decreasing by condition 3) 
$\n_1\ge \n_2 \ge \cdots \ge \n_j$ summing to $n-k$, for some $j$ with $1\le j\le k$ ($j$ depending on the number of singleton subtrees), and $\n_j\ge 2$ by condition 1. 
If gray edges from the root \emph{are} present,
then the root has no leaf child by conditions 2, 3, and 4, the gray edges are right justified, and the list of sizes of all subtrees of the root forms a partition $\n_1\ge \n_2 \ge \cdots \ge \n_j\ge 2$ of $n-k$ with $j>k$.
Putting all this together, we obtain the recurrence for $n\ge 2$,
\begin{equation}\label{eq1}
a_n =1 + \sum_{k=2}^{n-2}\,\sum_{(\n_1,\n_2,\dots,\n_j)\in \p_{n-k}^*} a_{\n_1}a_{\n_2}\cdots a_{\n_j}\, ,
\end{equation}
where $\p_r^*$ is the set $\{(\n_1,\n_2,\dots,\n_j): j\ge 1,\ \n_1\ge \n_2 \ge \cdots \ge \n_j \ge 2,\ \n_1+\n_2+\cdots +\n_j=r\}$ of partitions of $r$ with all parts $\ge 2$.
Replacing $n$ by $n+2$, (\ref{eq1}) becomes
\begin{equation}\label{eq2}
a_{n+2} =1 + \sum_{k=2}^{n}\,\sum_{(\n_1,\n_2,\dots,\n_j)\in \p_{n+2-k}^*} a_{\n_1}a_{\n_2}\cdots a_{\n_j}\, .
\end{equation}
In (\ref{eq2}), the summation index set $\mbf{\cup}_{k=2}^n \p_{n+2-k}^*$ consists of all partitions of integers between 2 and $n$ inclusive with all parts $\ge 2$. By appending 1's to these partitions (if necessary) we get all partitions of $n$ except the all-1s partition. Thus, by setting $a_1=1$, (\ref{eq2}) becomes, adopting the frequency-of-parts notation for integer partitions,
\begin{equation}\label{eq3}
a_{n+2}=\sum_{1^{i_1}2^{i_2}\dots\, n^{i_n}\vdash n}a_1^{i_1}a_2^{i_2} \cdots a_n^{i_n}\, ,
\end{equation} 
valid for $n\ge 1$, where the all-1s partition in (\ref{eq3}) produces the initial 1 on the right side of (\ref{eq2}). 
This is precisely the recurrence obtained by equating coefficients of $x^{n+2}$
in the defining generating function identity 
\begin{equation}\label{eq4}
\sum_{n\ge 1}a_n x^n = x+\frac{x^2}{\prod_{n\ge 1}(1-a_nx^n)}
\end{equation}
for \htmladdnormallink{A345973}{http://oeis.org/A345973} and (\ref{eq4}) immediately implies $a_1=1$.

\end{proof}
The trees in $\a_n$ with no gray edges are the trees counted by \htmladdnormallink{A346787}{http://oeis.org/A346787} (offset because A346787 counts by number of vertices which of course differs by 1 from the number of edges).

Department of Statistics, University of Wisconsin-Madison

\end{document}